\documentclass{lncse}
\def\Gr{Gr\"obner }
\begin{document}
\title{Constrained Hamiltonian Systems and \Gr
      Bases\thanks{
This work
was supported in part by Russian Foundation for Basic Research,
grant No. 98-01-00101.
} }
\author{Vladimir P. Gerdt\inst{1} \and Soso A. Gogilidze\inst{2}}

\institute{Laboratory of Computing Techniques and Automation,
Joint Institute for Nuclear Research,
141980 Dubna, Russia
\and
       Institute of High Energy Physics,
       Tbilisi State University,
       38086 Tbilisi, Georgia}
\maketitle
\begin{abstract}
In this paper we consider finite-dimensional constrained Hamiltonian
systems of polynomial type. In order to compute the complete set
of
constraints and separate them into the first and second classes we
apply the modern algorithmic methods of commutative algebra based on
the use of \Gr bases. As it is shown, this makes the classical Dirac
method fully algorithmic. The underlying algorithm implemented in Maple
is presented and some illustrative examples are given.
\end{abstract}

\section{Introduction}
The generalized Hamiltonian formalism invented by Dirac~\cite{Dirac}
for constrained systems has become a classical tool for investigation
of gauge theories in physics~\cite{GT,HT,PS}, and a platform for numerical
analysis of constrained mechanical systems~\cite{Seiler1}.
Finite-dimensional constrained Hamiltonian systems are part of
differential algebraic equations whose numerical analysis is of
great research interest over last decade~\cite{Campbell}
because of importance for many applied areas, for instance,
multi-body mechanics and molecular dynamics.

In physics, the constrained systems are mainly of interest for
purposes of quantization of gauge theories which play a fundamental role
in modern
quantum field theory and elementary particle physics. Dirac devised
his methods to study constrained Hamiltonian systems just for those
quantization purposes. Having this in mind, he classified the
constraints in the first and second classes. A first class constrained
physical system possesses gauge invariance and its quantization requires
gauge fixing whereas a second class constrained system does not need this.
The effect of the second class constraints may be reduced to a modification
of a naive measure in the path integral. The presence of gauge degrees of
freedom (first class constraints)
indicates that the general solution of the system depends on arbitrary
functions. Hence, the system is underdetermined. To eliminate unphysical
gauge degrees of freedom one usually imposes gauge fixing conditions
whereas for elimination of other unphysical degrees of freedom occurring
because of the second class constraints, one can use the Dirac
brackets~\cite{GT,HT,Sundermeyer}. In some special cases one can
explicitly eliminate the unphysical degrees of freedom~\cite{Soso1}.

\noindent
Unlike physics, where constrained systems are singular, as they
contain internal constraints, mechanical systems are usually
regular with externally imposed constraints~\cite{Arnold}.
Such a system is equivalent to a singular one whose Lagrangian is
that of the regular system enlarged with a linear combination of the
externally imposed constraints whose coefficients (multipliers) are
to be treated as extra dynamical variables. The latter system
may reveal extra constraints for the former system providing the consistency
of its dynamics.

Therefore, to investigate a constraint Hamiltonian system one has to detect
all the constraints involved, and separate them, for physical models,
into first and second classes. In his theory~\cite{Dirac} Dirac gave
the receipt for computation of constraints which is widely known as
{\em Dirac algorithm}, and it has been implemented in computer algebra
software~\cite{Alain}. However, the Dirac approach, as a method for
computation of constraints, is not yet
an algorithm. Even computation of the primary constraints, given a singular
Lagrangian, is not generally algorithmic. Moreover, in generation of the
secondary, tertiary, etc., constraints by the Dirac method one must verify
if a certain function of the phase space variables vanishes on the
constraint manifold. Generally, the latter problem is algorithmically
unsolvable. Similarly, there are no general algorithmic schemes for
separation of
constraints into the first and second classes.  In physical literature
one can find quite a number of particular methods developed for the
constraint separation (see, for example,~\cite{Lusanna,Pons}). But
all of them have non-algorithmic defects. Thereby, being successfully
applied to one constrained system, those methods may be failed for
another system even of a similar type.

In practice, many constrained physical and mechanical problems are
described by polynomial Lagrangians that lead to polynomial
Hamiltonians. In this case, as we show in the present paper, one can
apply \Gr bases which nowadays have become the most universal
algorithmic tool in commutative algebra~\cite{BW} and algebraic
geometry~\cite{CLO1,CLO2}. The combination of the Dirac method with
the \Gr bases technique makes the former fully algorithmic and, thereby,
allows to compute the complete set of constraints. Moreover, the
constraint separation is also done algorithmically. We show this
and present the underlying algorithm which we call {\em algorithm
Dirac-Gr\"{o}bner}. This algorithm has been implemented in
Maple V Release 5, and we illustrate it by examples both from physics
and mechanics.

\section{Dirac Method}

In this section we shortly describe the computational aspects of
the Dirac approach to constrained finite-dimensional Hamiltonian
systems~\cite{Dirac,HT}.

\noindent
Let us start with a Lagrangian
$L(q,\dot{q})\equiv L(q_i,\dot{q}_j)$ $(1\leq i,j\leq n)$
as a function of the generalized coordinates $q_i$ and velocities
$\dot{q}_j$\footnote{We consider only autonomous systems, and there is
no loss of generality since time $t$ may be treated as an additional
variable.}.
If the Hessian $\partial^2 L/{\partial \dot{q}_i}{\partial \dot{q}_j}$
has the full rank $r=n$, then the system
is {\em regular} and it has no internally hidden constraints.
Otherwise, if $r<n$,
the Euler-Lagrange equations
\begin{equation}
\dot{p}_i=\frac{\partial L}{\partial q_i}\quad (1\leq i\leq n) \label{Lag_eq}
\end{equation}
with
\begin{equation}
p_i=\frac{\partial L}{\partial \dot{q}_i} \label{def_p}
\end{equation}

\noindent
are {\em singular} or {\em degenerate}, as
not all differential equations~(\ref{Lag_eq}) are of the second order.
There are
just $n-r$ such independent lower order equations.
By the Legendre
transformation~\footnote{In this paper summation
over repeated indices is we always assumed.}
\begin{equation}
H_c(p,q)=p_iq_i - L \label{def_H_c},
\end{equation}
we obtain the {\em canonical Hamiltonian} with momenta $p_i$
defined in~(\ref{def_p}). In the degenerate case there are
{\em primary constraints}
denoted by $\phi_\alpha$, which form the {\em primary constraint manifold}
denoted by $\Sigma_0$
\begin{equation}
\Sigma_0\ :\quad \phi_\alpha (p,q)=0\quad  (1\leq \alpha\leq n-r),
\label{pr_constr}
\end{equation}
Thus, the dynamics of the system is determined only on the
constraint manifold~(\ref{pr_constr}).  To take this fact into account,
Dirac defined the {\em total Hamiltonian}
\begin{equation}
 H_t=H_c+u_\alpha \phi_\alpha  \label{def_h_t}
\end{equation}
with {\em multipliers} $u_\alpha$ as arbitrary (non-specified) functions
of the coordinates and momenta. The corresponding Hamiltonian equations
determine the system dynamics together with the primary constraints
\begin{equation}
\dot{q}_i=\{H_t,q_i\},\ \ \dot{p}_i=\{H_t,p_i\},\ \
\phi_\alpha(p,q)=0\ \ (1\leq i\leq n,\ 1\leq \alpha\leq n-r), \label{H_eq}
\end{equation}
where the {\em Poisson brackets} are defined for any two functions $f,g$
of the dynamical variables $p$ and $q$ as follows
\begin{equation}
\{f,g\}=\frac{\partial f}{\partial p_i} \frac{\partial g}{\partial q_i}-
\frac{\partial g}{\partial p_i} \frac{\partial p}{\partial q_i}. \label{p_b}
\end{equation}
In order to be consistent with the system dynamics, the primary constraints
must satisfy the conditions
\begin{equation}
\dot{\phi}_\alpha=\{H_t,\phi_\alpha\} \stackrel{\Sigma_0}{=} 0\quad
(1\leq \alpha\leq n-r), \label{cons_cond}
\end{equation}

\noindent
where $\stackrel{\Sigma_0}{=}$ stands for the equality, called
{\em a week equality}, on the primary
constraint manifold~(\ref{pr_constr}). The
Poisson bracket in (\ref{cons_cond}) must be a linear combination of
the constraint functions~\cite{HT}. Given a constraint function
$\phi_\alpha$,
the consistency condition~(\ref{cons_cond}), unless it is satisfied
identically, may lead either to a contradiction or to a new constraint.
The former case signals that the given Hamiltonian system is inconsistent.
In the latter case, if the new constraint does not involve any of
multipliers
$u_\alpha$, it must be added to the constraint set, and, hence, the
constraint manifold must involve this new constraint. Otherwise, the
consistency condition is considered as defining the multipliers, and the
constraint set is not enlarged with it.

The iteration of this consistency check ends up with the {\em complete
set of constraints} such that for every constraint in the set condition
(\ref{cons_cond}) is satisfied. This is the Dirac method of the constraint
computation. As shown in~\cite{Seiler2}, the method is nothing
else than completion of the initial Hamiltonian system to involution,
and the
constraints generated are just {\em the integrability conditions}. For
general systems of PDEs, the completion process is done~\cite{Pommaret}
by sequential prolongations and projections. For Hamiltonian systems,
the time derivative of a constraint
is its prolongation whereas projection of the prolonged constraint
is realized in (\ref{cons_cond}) by computing the Poisson bracket on the
constraint manifold.

Let now $\Sigma$ be the constraint manifold for the complete set of
constraints
\begin{equation}
\Sigma\ :\quad \phi_\alpha (p,q)=0\quad  (1\leq \alpha\leq k).
\label{constr_man}
\end{equation}

\noindent
If a constraint function $\phi_\alpha$ satisfies the condition
\begin{equation}
\{\phi_\alpha(p,q),\phi_\beta(p,q)\} \stackrel{\Sigma}{=} 0\quad
(1\leq \beta\leq k), \label{fc_cons}
\end{equation}

\noindent
it is of {\em the first class}. Otherwise, the constraint
function is of {\em the second class}. The number of the second class
constrains is equal to rank of the following
$(k\times k)$ {\em Poisson bracket matrix}, whose elements must be
evaluated on the constraint manifold

\begin{equation}
M_{\alpha \beta} \stackrel{\Sigma}{=} \{\phi_\alpha,\phi_\beta \}.
\label{PBM}
\end{equation}
Note that matrix $M$ has even rank because of its skew-symmetry.

If a Lagrangian system $L_0(q,\dot{q})$ is regular with externally
imposed {\em holonomic} constraints $\psi_\alpha(q)=0$, the system is
equivalent~\cite{Seiler1} to the singular one with Lagrangian
$L=L_0+\lambda_\alpha \phi_\alpha$ and extra generalized
coordinates $\lambda_\alpha$. Furthermore, the Dirac method can be applied
for finding the other constraints inherent in the initial regular system
and, hence, not involving the extra dynamical variables.

Therefore, the problem of constraint computation and separation is
reduced
to manipulation with functions of the coordinates and momenta on the
constraint manifold. Generally, there is no algorithmic way for such
a manipulation. However, for polynomial functions all the related
computations can be done algorithmically by means of \Gr bases,
as we show in the next section.

\section{Algorithm Description}
Here we describe an algorithm which,
given a polynomial Lagrangian whose coefficients are rational numbers,
computes the complete set of constraints and separates them into
the first and second classes. The algorithm combines the above described
Dirac method with the \Gr bases technique. By this reason we call it
algorithm Dirac-Gr\"{o}bner. All the below used concepts, definitions and
constructive methods related to \Gr bases are explained, for instance,
in textbooks \cite{BW,CLO1,CLO2}.

At first we present the algorithm under assumption that a polynomial
ideal generated by constraints is radical. This is true
for most of real practical problems. Next, we indicate
how to modify the algorithm to treat the most general (non-radical)
case.
\vskip 0.3cm
\centerline{\bf Algorithm Dirac-\Gr }
\vskip 0.2cm
\noindent
{\bf Input:} $L(q,\dot{q})$, a polynomial Lagrangian $(L\in Q[q,\dot{q}])$
\vskip 0.1cm
\noindent
{\bf Output:} $\Phi_1$ and $\Phi_2$, sets of the first and second class
constraints, respectively.
\begin{enumerate}
\item Computation of the canonical Hamiltonian
and primary constraints:
\vskip 0.1cm
\begin{enumerate}
\item Construct the polynomial set $F=\cup_{i=1}^n
 \{p_i-\partial L/\partial \dot{q}_i\}$ in variables $p,q,\dot{q}$.
\item Compute the \Gr basis $G$ of the ideal in
ring $Q[p,q,\dot{q}]$ generated by $F$ with respect to
an ordering\footnote{An elimination ordering which induced the
degree-reverse-lexicographical one for monomials in $p$ and $q$ is
heuristically best for efficiency reasons.} which eliminates $\dot{q}$.
Then compute
the canonical Hamiltonian as the normal form of (\ref{def_H_c}) modulo $G$.
\item Find the set $\Phi$ of primary constraint polynomials
as $G\cap Q[p,q]$.  If $\Phi=\emptyset$, then stop since the system is
regular. Otherwise, go to the next step.
\end{enumerate}

\vskip 0.2cm
\item Computation of the complete set of constraints:

\vskip 0.1cm
\begin{enumerate}
\item Take $G=\Phi$ for the \Gr basis $G$ of the ideal generated by $\Phi$ in
$Q[p,q]$ with respect to the ordering induced by that chosen at
Step 1(b). Fix this ordering in the sequel.
\item Construct the total Hamiltonian in form (\ref{def_h_t}) with
multipliers $u_\alpha$
treated as symbolic constants (parameters).
\item For every element $\phi_\alpha$ in $\Phi$ compute the
normal form $h$ of the Poisson bracket $\{H_t,\phi_\alpha\}$
modulo $G$. If $h\neq 0$ and no multipliers $u_\beta$ occur in it, then
enlarge set $\Phi$ with $h$, and compute the \Gr basis $G$ for the
enlarged set.
\item If $G=\{1\}$, stop because the system is inconsistent. Otherwise,
repeat the previous step until the consistency condition (\ref{cons_cond})
is satisfied for every element in $\Phi$ irrespective of
multipliers $u_\alpha$. This gives the complete set of constraints
$\Phi=\{\phi_1,\ldots,\phi_k\}$.
\end{enumerate}

\vskip 0.2cm
\item Separation of constraints into first and second
classes:

\vskip 0.1cm
\begin{enumerate}
\item Construct matrix $M$ in (\ref{PBM}) by computing the
normal forms of its elements modulo $G$, and determine
rank $r$ of $M$. If $r=k$, stop with $\Phi_1=\emptyset$, $\Phi_2=\Phi$. If
$r=0$, stop with $\Phi_1=\Phi$ and $\Phi_2=\emptyset$. Otherwise,
go to the next step.
\item Find a basis $A=\{a_1,\ldots,a_{k-r}\}$ of the null space (kernel)
of the linear transformation defined by $M$. For every vector $a$ in $A$
construct a first class
constraint as $a_\alpha \phi_\alpha$. Collect them
in set $\Phi_1$.
\item Construct $(k-r)\times k$ matrix $(a_j)_\alpha$ from components of
vectors in $A$ and find a basis $B=\{b_1,\ldots,b_{r}\}$ of the
null space of the corresponding linear transformation. For every vector
$b$ in $B$ construct a second
class constraint as $b_\alpha \phi_\alpha$. Collect them
in set $\Phi_2$.
\end{enumerate}
\end{enumerate}
The correctness of Steps 1, 2 and 3(a) of the algorithm is provided
by the properties of \Gr bases~\cite{BW,CLO1,CLO2} and by the following
facts:
(i) the definition (\ref{def_H_c}) of the canonical Hamiltonian implies
its independence of $\dot{q}$ on the primary constraint manifold
(\ref{pr_constr}); (ii) whenever a multiplier $u_\alpha$ in
(\ref{def_h_t}) is differentiated when the Poisson
bracket in (\ref{cons_cond}) is evaluated, the corresponding term
vanishes on the constraint manifold. The correctness of Steps 3(b)
and 3(c) follows
from definition (\ref{fc_cons}) of the first class constraints and
the correctness of Step 3(a). The termination of algorithm
Dirac-Gr\"{o}bner follows from the finiteness of the \Gr basis $G$
which is constructed at Step 2(c).

Now consider the most general case when the constraints obtained
from (\ref{cons_cond}) lead to a non-radical ideal. It should be noted
that the ideal generated by the primary constraint polynomials (Step 1)
is always radical. This is provided by linearity of (\ref{def_p}) in momenta.
However, already the first secondary constraint added may destroy this
property of the ideal. Therefore, the algorithm needs one more step, namely,
Step 2(e), where the \Gr basis $G$ of the radical ideal for the polynomial
set
$\Phi$ is computed. Next, every constraint polynomial in
$\Phi$ is replaced by its normal form modulo $G$. All the elements
with zero normal forms are eliminated from the set. The extra step is
also algorithmic. There are algorithms for construction of a basis,
and, hence, a \Gr basis, of the radical of a given ideal, which are
built-in in some computer algebra systems (see~\cite{BW,CLO1,CLO2} for
more details and references). One can also check the radical
membership of $h$ at Step 2(c) before its adding to $\Phi$. This check
is easily done~\cite{BW,CLO1}, but in any case Step 2, for the correctness
of Step 3, must end up with the radical sets $\Phi$ and $G$.

We implemented algorithm Dirac-Gr\"{o}bner, as it presented above
for the radical case, in Maple V Release 5. The implementation is relied
on the built-in system facilities for computation and manipulation with
\Gr bases and for linear algebra. Using our Maple code for different
examples from physics and mechanics, we experimentally observed that
in those infrequent cases when the constraint ideals are non-radical
this can easily be detected from the structure of the output set.

\section{Examples}

In this section we illustrate, by examples from physics and mechanics,
the application of algorithm Dirac-Gr\"{o}bner.

\begin{example} $SU(2)$ Yang-Mills mechanics in $0+1$ dimensional
space-time~\cite{Soso1}. This is a constrained physical
model with gauge symmetry. The model Lagrangian is given by
$L=\frac{1}{2}(D_t)_i(D_t)_i$, $(D_tx)_i=\dot{x}_i+g\epsilon_{ijk}y_jx_k$
$(1\leq i,j,k\leq 3)$.
Here $x_i$ and $y_i$ are the generalized coordinates and tensor
$\epsilon_{ijk}$ is anti-symmetric in its indices with
$\epsilon_{123}=1$. Respectively, the primary constraints and the
canonical Hamiltonian are $p_i^y=0$ and
$H_c=\frac{1}{2}-\epsilon_{ijk}x_jp_ky_i$ with the momenta given
by $p_i^y=\partial L/\partial \dot{y}_i$ and
$p_i=\partial L/\partial \dot{x}_i$. The other constraints in the
complete set computed by the algorithm are $\phi_i=\epsilon_{ijk}x_jp_k=0$,
and all the six constraints found are of the first class.
\end{example}

\begin{example} Point particle of mass $m$ moving on the surface of
a sphere (rigid rotator). The movement is described
by the regular Lagrangian
$L_0=\frac{1}{2}m^2(\dot{q_1}^2+\dot{q_2}^2+\dot{q_3}^2)/2
\equiv \frac{1}{2}m^2\dot{q}^2$
with the externally imposed holonomic constraint
$\phi(q)=q^2-1=0$.
This system is equivalent to the singular Lagrangian system
$ L=L_0 + \lambda \phi$, where $\lambda$ is an extra coordinate.
There is the only primary constraint $p_\lambda=0$
$(p_\lambda=\partial L/\partial \lambda)$, and the canonical
Hamiltonian is $H_c=\frac{1}{2}m^2p^2-\lambda \phi(q)$
$(p_i=\partial L/\partial q_i)$. The complete set of constraint
polynomials for the singular system contains four second class polynomials
$\{p_\lambda,\phi(q),p_iq_i,2m\lambda + p^2\}$. Coming back to the initial
regular system, the first and the last polynomials in the set must be
omitted since they determine the extra dynamical variables.
\end{example}

\begin{example} Singular physical system with both first and second
class constraints\footnote{A.Burnel. Private communication.}.
The system Lagrangian is $L=q_1(\dot{q}_2-q_3) - \dot{q}_1q_2$.
There are three primary constraint polynomials $\{p_1+q_2,p_2-q_1,p_3\}$.
The canonical Hamiltonian is $H_c=q_1q_2$. One more constraint polynomial
$q_1$ is found by the
Dirac-Gr\"{o}bner algorithm. The sets $\Phi_1$ and $\Phi_2$ of the
first and second classes are $\{p_2+q_1,p_3\}$ and $\{p_1+q_2,q_1\}$,
respectively. Note that this system has no physical degrees of
freedom (c.f.~\cite{Seiler2}).
\end{example}

\begin{example} Inconsistent singular
system~\cite{PS}:  $L=\frac{1}{2}\dot{q}_1^2+q_2$.
There is the single primary constraint $p_2=0$. The canonical
Hamiltonian is $H_c=p_1^2/2-q_2$. At Step 2(c) of
algorithm Dirac-Gr\"{o}bner the inconsistency $\dot{p}_2=1$ occurs.
The algorithm detects this inconsistency and stops.
\end{example}

\noindent
The above examples are rather small and can be treated by hand. With our
Maple code we have already tried successfully much more nontrivial
examples. For instance, we computed and separated the constraints for
the $SU(2)$ Yang-Mills mechanics in $3+1$ dimensional
space-time~\cite{Soso1}. Surprisingly, this computation took
only a few seconds on an Pentium 100 personal computer though
the model Lagrangian and the canonical Hamiltonian are rather
cumbersome polynomials of the 4th degree in 21 variables.

\end{document}